\documentclass[11pt,reqno]{amsart}

\usepackage{mathtools,amssymb,amsfonts,amsthm,booktabs,enumitem}
\usepackage{xcolor}
\usepackage[colorlinks,linkcolor=magenta,anchorcolor=blue,citecolor=green]{hyperref}

\newtheorem{theorem}{Theorem}[section]
\newtheorem{lemma}[theorem]{Lemma}
\newtheorem{proposition}[theorem]{Proposition}
\newtheorem{corollary}[theorem]{Corollary}
\theoremstyle{definition}
\newtheorem{question}[theorem]{Question}
\newtheorem{hypothesis}[theorem]{Hypothesis}
\newtheorem{conjecture}[theorem]{Conjecture}
\newtheorem{construction}[theorem]{Construction}
\newtheorem{remark}[theorem]{Remark}

\newcommand\Cay{\mathrm{Cay}}
\newcommand\Haar{\mathrm{Haar}}
\newcommand\Dih{\mathrm{Dih}}
\newcommand\D{D}
\newcommand\Z{Z}
\newcommand\ind{\mathbf{1}}

\begin{document}

\title[Hamilton cycles in generalized dihedral Cayley digraphs]{Hamilton cycles in generalized dihedral Cayley graphs and digraphs}

\author{Jiyong Chen}
\address{School of Mathematical Sciences, Xiamen University, Xiamen 361005, China}
\email{chenjy1988@xmu.edu.cn}

\author{Jianyu Ou}
\address{School of Mathematical Sciences, Xiamen University, Xiamen 361005, China}
\email{oujianyu@xmu.edu.cn}

\author{Yan-Li Qin}
\address{School of Statistics, Capital University of Economics and Business, Beijing 100070, China}
\email{ylqin@cueb.edu.cn}

\author{Binzhou Xia}
\address{School of Mathematics and Statistics, The University of Melbourne, Parkville, VIC 3010, Australia}
\email{binzhoux@unimelb.edu.au}

\author{Kai Yuan}
\address{School of Mathematics and Information Science, Yantai University, Yantai 264005, China}
\email{pktide@163.com}

\begin{abstract}
We prove that every connected Cayley digraph on a generalized dihedral group of order at least $4$ has a directed Hamilton cycle. In particular, this confirms a conjecture of Holszty\'{n}ski and Strube from 1978 for dihedral groups. The key new ingredient is a three-fold sumset covering theorem for the terminal coordinates of Hamilton paths in cubic Haar graphs over abelian groups of odd order, with connection sets minimal subject to connectivity.

\textit{Key words:} Hamilton cycle; Cayley graph; Cayley digraph; generalized dihedral group; Haar graph.

\textit{MSC2020:} 05C45, 05C25, 05C20.
\end{abstract}

\maketitle

\section{Introduction}\label{sec:intro}

All groups and (di)graphs in this paper are assumed to be finite unless otherwise stated. For a group $G$ and a subset $S\subseteq G\setminus\{1\}$, the \emph{Cayley digraph} $\Cay(G,S)$ has vertex set $G$ and arcs $(g,sg)$ for $g\in G$ and $s\in S$. When $S$ is inverse-closed, $\Cay(G,S)$ is a \emph{Cayley graph}, where we identify each pair of opposite arcs with an undirected edge. Recall that a digraph is called connected if its underlying undirected graph is connected. For $\Cay(G,S)$, this is equivalent to $\langle S\rangle=G$.

Lov\'{a}sz's question from 1969~\cite{Lovasz1969}, whether every connected vertex-transitive graph has a Hamilton path, led to the following famous question in the study of Cayley graphs.

\begin{question}\label{que:cayley}
Does every connected Cayley graph with more than two vertices have a Hamilton cycle?
\end{question}

The answer is affirmative for abelian groups. More precisely, Chen and Quimpo~\cite{CQ1980} proved that every connected Cayley graph of valency at least $3$ on an abelian group is Hamilton-connected or Hamilton-laceable. Here a graph is called \emph{Hamilton-connected} if every two distinct vertices are joined by a Hamilton path, and a bipartite graph is called \emph{Hamilton-laceable} if every two vertices in different parts are joined by a Hamilton path. Either property gives a Hamilton cycle when the graph has more than two vertices.

Positive answers are also known for semidirect products of a cyclic group of prime order by an abelian group~\cite{Durnberger1983,Marusic1983}, for nilpotent groups with a cyclic commutator subgroup~\cite{GM2014}, and for groups whose commutator subgroup is cyclic of prime-power order~\cite{KW1985}. Further results are surveyed in~\cite{CG1996,KM2009,WG1984}. For dihedral groups, Holszty\'{n}ski and Strube~\cite[page~270]{HS1978} proposed the following conjecture in 1978.

\begin{conjecture}\label{conj:dihedral}
Every connected Cayley digraph on a dihedral group has a directed Hamilton cycle.
\end{conjecture}

Holszty\'{n}ski and Strube~\cite[Proposition~6.1]{HS1978} proved the conjecture for twice prime order. Alspach and Zhang~\cite{AZ1989} settled the undirected cubic case. Further progress applies to generalized dihedral groups. Here, the \emph{generalized dihedral group} $\Dih(A)$ of a nontrivial abelian group $A$ is the semidirect product of $A$ by $\Z_2$, where the nonidentity element of $\Z_2$ inverts every element of $A$. In particular, $\Dih(A)$ is dihedral when $A$ is cyclic.

Alspach, Chen and Dean~\cite[Theorem~1.8]{ACD2010} proved that every connected Cayley graph of valency at least $3$ on a generalized dihedral group of order divisible by $4$ is Hamilton-connected or Hamilton-laceable. Pastine and Jaume~\cite[Theorem~15]{PJ2012} proved that the connected Cayley digraph $\Cay(\Dih(A),S)$ has a directed Hamilton cycle whenever $S\cap A\neq\varnothing$. Since every element of $\Dih(A)\setminus A$ is an involution and a connected graph of valency $2$ is a cycle, these results leave only the case where $|A|$ is odd and $S\subseteq\Dih(A)\setminus A$.

In this paper we settle the remaining case and prove the following theorem.

\begin{theorem}\label{thm:generalized}
Every connected Cayley digraph on a generalized dihedral group of order at least $4$ has a directed Hamilton cycle.
\end{theorem}

In particular, Theorem~\ref{thm:generalized} answers Question~\ref{que:cayley} affirmatively for generalized dihedral groups, since Cayley graphs are a special case of Cayley digraphs. It also confirms the 1978 conjecture of Holszty\'{n}ski and Strube.

\begin{corollary}\label{cor:dihedral}
Conjecture~\ref{conj:dihedral} holds.
\end{corollary}

\section{Preliminaries}\label{sec:prelim}

Throughout the paper, $\log$ denotes the natural logarithm (to base $e$). For a property $\mathcal P$, write $\ind_{\mathcal P}=1$ if $\mathcal P$ holds, and $\ind_{\mathcal P}=0$ otherwise.
For an integer $n\geq1$, let $\omega(n)$ be the number of distinct prime divisors of $n$. Let $\phi$ be Euler's totient function, and let $\mu$ be the M\"{o}bius function, which is defined by
\[
\mu(n)=
\begin{cases}
0&\text{if }n\text{ is divisible by the square of a prime}\\
(-1)^r&\text{if }n\text{ is a product of }r\text{ distinct primes}.
\end{cases}
\]
In particular, $\omega(1)=0$ and $\phi(1)=\mu(1)=1$. For an integer $k\geq1$, we have the basic identity
\begin{equation}\label{eq:mobius-divisor-sum}
\sum_{d\mid k}\mu(d)=
\begin{cases}
1&\text{if }k=1\\
0&\text{if }k>1,
\end{cases}
\end{equation}
where the sum runs over all the positive divisors $d$ of $k$. Let
\[
\Phi(n)=\sum_{i=1}^n\phi(i).
\]
Although it is well known that $\Phi(n)=3n^2/\pi^2+O(n\log(n))$ as $n\to\infty$ (see for example~\cite[Theorem~3.7]{Apostol1976}), we need the following lower bound for $\Phi(n)$ that is valid for every $n\geq1$.

\begin{lemma}\label{lem:totient}
$\Phi(n)\geq2n(n+1)/7$ for every integer $n\geq1$.
\end{lemma}

\begin{proof}
Exact computation verifies the inequality for $n\leq1000$ (with equality precisely at $n=6$). Suppose that $n\geq1000$.  By~\eqref{eq:mobius-divisor-sum}, for positive integers $a$ and $b$,
\begin{equation}\label{eq:indicator}
\sum_{d\mid\gcd(a,b)}\mu(d)=\ind_{\gcd(a,b)=1}.
\end{equation}
We count the ordered coprime pairs in $\{1,\ldots,n\}\times\{1,\ldots,n\}$. There are
\[
\sum_{b=1}^n\phi(b)=\Phi(n)
\]
such pairs with $a\leq b$. Hence the total is $2\Phi(n)-1$, as $(a,b)=(1,1)$ is the only coprime pair with $a=b$. Then applying~\eqref{eq:indicator} and interchanging the sums gives
\begin{equation}\label{eq:Phi}
2\Phi(n)-1
=\sum_{a=1}^n\sum_{b=1}^n\sum_{d\mid\gcd(a,b)}\mu(d)
=\sum_{d=1}^n\mu(d)\left\lfloor\frac{n}{d}\right\rfloor^2
>n^2\sum_{d=1}^n\frac{\mu(d)}{d^2}-2n\sum_{d=1}^n\frac1d.
\end{equation}
The second equality counts the pairs for which $d$ divides both coordinates. For the inequality, use $\lfloor x\rfloor^2>x^2-2x$ when $\mu(d)=1$ and $-\lfloor x\rfloor^2\geq-x^2$ when $\mu(d)=-1$.

Recall that the Dirichlet series identity for the M\"{o}bius function (see~\cite[Section~11.4]{Apostol1976}) gives
\[
\sum_{d=1}^{\infty}\frac{\mu(d)}{d^2}=\frac{1}{\zeta(2)}.
\]
As a consequence,
\begin{align*}
\sum_{d=1}^n\frac{\mu(d)}{d^2}
=\frac1{\zeta(2)}-\sum_{d=n+1}^{\infty}\frac{\mu(d)}{d^2}
\geq\frac1{\zeta(2)}-\sum_{d=n+1}^{\infty}\frac1{d^2}
>\frac1{\zeta(2)}-\int_n^{\infty}\frac{dx}{x^2}
=\frac6{\pi^2}-\frac1n.
\end{align*}
Combining this with~\eqref{eq:Phi} and $\sum_{d=1}^n1/d<1+\log(n)$, we conclude that
\[
\Phi(n)>\frac{3n^2}{\pi^2}-\frac n2-n\sum_{d=1}^n\frac1d+\frac12
>\frac{3n^2}{\pi^2}-n\Big(\log(n)+\frac32\Big)+\frac12,
\]
where the last lower bound exceeds $2n(n+1)/7$ for $n\geq1000$. Indeed, let $c=3/\pi^2-2/7$. The difference is
\[
n\left(cn-\log(n)-\frac{25}{14}\right)+\frac12.
\]
Since $\pi<22/7$, we have $c>61/3388>1/1000$, so $cn-\log(n)-25/14$ is increasing for $n\geq1000$. Its value at $1000$ is positive, as
\[
1000c-\log(1000)-\frac{25}{14}
>\frac{61000}{3388}-7-\frac{25}{14}>0.
\]
This proves the required bound.
\end{proof}

For a subset $S$ of a group $G$, the \emph{Haar graph} $\Haar(G,S)$ has vertex set $G\times\{0,1\}$ and edge set
$\{\{(g,0),(h,1)\}\mid g,h\in G,\,hg^{-1}\in S\}$. Write $g_+=(g,0)$ and $g_-=(g,1)$ for $g\in G$, and call $g$ the \emph{coordinate} of the vertices $g_+$ and $g_-$. An edge of $\Haar(G,S)$ has the form
\[
\{g_+,(sg)_-\}
\]
for some $g\in G$ and $s\in S$, and is called an \emph{$s$-edge}. Write $U_\epsilon=\{g_\epsilon\mid g\in U\}$ for $U\subseteq G$ and $\epsilon\in\{+,-\}$. Then $\Haar(G,S)$ is bipartite with parts $G_+$ and $G_-$, and
\[
g_+\sim h_-\ \ \iff\ \ hg^{-1}\in S.
\]
When $G$ is abelian, we use additive notation. We write $\Haar(n,S)$ for $\Haar(\Z_n,S)$, which are the \emph{cyclic Haar graphs} studied in~\cite{HMP2002}.

For a subset $E$ of an abelian group and a positive integer $m$, write
\[
mE=E+\cdots+E=\{x_1+\cdots+x_m\mid x_1,\ldots,x_m\in E\}
\]
for the $m$-fold sumset, allowing repeated summands. For an abelian group $A$, let
\[
\Dih(A)=A\rtimes\langle\sigma\rangle,
\]
where $\sigma$ is an involution acting on $A$ by inversion. In particular,
\[
\Dih(Z_n)=\D_{2n}=\langle\rho,\sigma\mid\rho^n=\sigma^2=1,\ \sigma\rho\sigma=\rho^{-1}\rangle.
\]

\begin{lemma}\label{lem:haar}
Let $A$ be an abelian group, and let $S\subseteq A$. The following statements hold:
\begin{enumerate}[label=\textup{(\alph*)},leftmargin=*]
\item\label{item:haar-cayley} $\Haar(A,S)$ is isomorphic to $\Cay(\Dih(A),S\sigma)$ via the map $a_{+}\mapsto-a,\,a_{-}\mapsto a\sigma$.
\item\label{item:haar-translation} $\Haar(A,S)$ is isomorphic to $\Haar(A,S+c)$ for each $c\in A$ via the map $a_{+}\mapsto a_{+}$, $a_{-}\mapsto (a+c)_{-}$.
\item\label{item:translation} For each $c\in A$, translation of $c$ induces an automorphism of $\Haar(A,S)$ via the map $a_{+}\mapsto(a+c)_+,\,a_{-}\mapsto(a+c)_-$.
\item\label{item:haar-components} Let $B=\langle s-t\mid s,t\in S\rangle$. For each $r\in S$, the component of $\Haar(A,S)$ containing $0_{+}$ has the biparts $B_+$ and $(B+r)_-$, and is isomorphic to $\Haar(B,S-r)$ via the map $b_{+}\mapsto b_{+},\,b_{-}\mapsto (b-r)_{-}$. If $S\neq\varnothing$, then $\Haar(A,S)$ has $|A/B|$ connected components.
\end{enumerate}
\end{lemma}

\begin{proof}
It is straightforward to verify~\ref{item:haar-cayley},~\ref{item:haar-translation} and~\ref{item:translation}.
For~\ref{item:haar-components}, let $s\in S$. A walk of length two from $a_{+}$ ends at $(a+s-t)_{+}$ for some $s,t\in S$, and every such step is possible. Hence the component of $\Haar(A,S)$ containing $0_{+}$ intersects $A_+$ at $B_+$, and the neighbors of vertices in $B_+$ form $(B+r)_-$. Moreover, the stated map is readily verified to be an isomorphism. Finally, by~\ref{item:translation}, the automorphisms of $\Haar(A,S)$ that are translations by $c\in A$ permute the components. This gives $|A/B|$ components if $S\neq\varnothing$.
\end{proof}

For a subset $S$ of an abelian group $A$, denote
\[
D(S)=\langle s-t\mid s,t\in S\rangle.
\]
By Lemma~\ref{lem:haar}\ref{item:haar-components}, $\Haar(A,S)$ is connected if and only if $S\neq\varnothing$ and $D(S)=A$. Denote
\[
\mathrm{HP}(A,S)=\{a\in A\mid \Haar(A,S)\text{ has a Hamilton path from }0_{+}\text{ to }a_{-}\},
\]
and write $\mathrm{HP}(n,S)=\mathrm{HP}(Z_n,S)$. For each $c\in A$, Lemma~\ref{lem:haar}\ref{item:haar-translation} implies
\[
\mathrm{HP}(A,S+c)=\mathrm{HP}(A,S)+c.
\]

In the study of Hamilton paths or cycles in a Haar graph, we usually need to join Hamilton paths in some smaller Haar graphs.

\begin{hypothesis}\label{hyp:concat}
Let $A=B+\langle s\rangle$ be a finite abelian group, where $B$ is a subgroup of index $m\geq2$ and $s\in A$, let $T\subseteq B$ be nonempty with $D(T)=B$, and let $S=T\cup\{s\}$. For each $i\in \{0,1,\ldots,m-1\}$, let $B_i=B-is$.
\end{hypothesis}

Under Hypothesis~\ref{hyp:concat}, Lemma~\ref{lem:haar}\ref{item:haar-components} shows that $(B_i)_+\cup(B_i)_-$ with $i\in \{0,1,\ldots,m-1\}$ are the $m$ components of $\Haar(A,T)$. We use the following construction of a Hamilton path in $\Haar(A,S)$ from Hamilton paths in $\Haar(B,T)$.

\begin{construction}\label{con:concat}
Under Hypothesis~\ref{hyp:concat}, for Hamilton paths $P_0,P_1,\ldots,P_{m-1}$ in $\Haar(B,T)$ from $0_{+}$ to $(b_0)_{-},(b_1)_-,\ldots,(b_{m-1})_{-}$, respectively, define $a_0,a_1,\ldots,a_{m-1}$ inductively by letting
\begin{equation}\label{eq:a}
a_0=0,\quad a_i=a_{i-1}+b_{i-1}-s
\end{equation}
for $i\in\{1,\ldots,m-1\}$. For each $i\in\{0,1,\ldots,m-1\}$, it follows that $B_i=B+a_i$. Moreover, let
\begin{equation}\label{eq:c}
c_i=a_i+b_i,
\end{equation}
and let $Q_i$ be the path obtained from $P_i$ by applying the translation of $a_i$. Then $Q_i$ is a Hamilton path in the component $(B_i)_+\cup(B_i)_-$ of $\Haar(A,T)$, which starts at $(a_i)_{+}$ and ends at $(c_i)_{-}$. For each $i\in\{1,\ldots,m-1\}$, we derive from~\eqref{eq:a} and~\eqref{eq:c} that $a_i=c_{i-1}-s$, and so $\{(c_{i-1})_{-},(a_i)_{+}\}$ is an $s$-edge of $\Haar(A,S)$ joining $Q_{i-1}$ to $Q_i$. The resulting concatenation $Q_0Q_1\cdots Q_{m-1}$ is a Hamilton path of $\Haar(A,S)$ from $0_{+}$ to $(c_{m-1})_{-}$, where
\begin{equation}\label{eq:delta}
c_{m-1}=\sum_{i=0}^{m-1}b_i-(m-1)s.
\end{equation}
The $s$-edges of $Q_0Q_1\cdots Q_{m-1}$ are precisely the $m-1$ joining edges $\{(c_{i-1})_{-},(a_i)_{+}\}$ with $i\in\{1,\ldots,m-1\}$.
\end{construction}

\begin{lemma}\label{lem:concat}
Under Hypothesis~\ref{hyp:concat}, suppose that $\ell\mathrm{HP}(B,T)=B$ for some positive integer $\ell\leq m$. The following statements hold:
\begin{enumerate}[label=\textup{(\alph*)},leftmargin=*]
\item\label{item:concat-coverage} $B+s\subseteq\mathrm{HP}(A,S)$.
\item\label{item:concat-cycle} $\Haar(A,S)$ has a Hamilton cycle containing $\{0_+,s_-\}$.
\end{enumerate}
\end{lemma}

\begin{proof}
\ref{item:concat-coverage} Since $\mathrm{HP}(B,T)$ is nonempty and $B+\mathrm{HP}(B,T)=B$, the assumption $\ell\mathrm{HP}(B,T)=B$ with $\ell\leq m$ implies $m\mathrm{HP}(B,T)=B$. Choose $P_0,P_1,\ldots,P_{m-1}$ freely in Construction~\ref{con:concat}. The sum of their terminal coordinates ranges over $B$. Hence~\eqref{eq:delta} shows that $c_{m-1}$ ranges over $B-(m-1)s=B+s$, as $ms\in B$. This proves the inclusion.

\ref{item:concat-cycle} By statement~\ref{item:concat-coverage}, there is a Hamilton path from $0_+$ to $s_-$. Since $|A|\geq2$, the path has at least four vertices, so the edge $\{0_+,s_-\}$ does not belong to it. Adding this edge gives the required Hamilton cycle.
\end{proof}

Let $P=w_0w_1\cdots w_\ell$ be a path in a graph, and suppose that $w_\ell$ is adjacent to $w_i$ for some $i\in\{0,1,\ldots,\ell-2\}$. Replacing the edge $\{w_i,w_{i+1}\}$ by $\{w_i,w_\ell\}$ gives the path
\[
w_0\cdots w_i\,w_\ell w_{\ell-1}\cdots w_{i+1}.
\]
This is the \emph{P\'{o}sa rotation}~\cite{Posa1976} of $P$ at the \emph{pivot} $w_i$. It preserves the initial vertex and the vertex set, and its terminal vertex is the successor of the pivot on $P$.

\section{Hamilton paths in cubic Haar graphs}\label{sec:three}

The main purpose of this section is to prove the following proposition.

\begin{proposition}\label{prop:threefold}
Let $A$ be a finite abelian group of odd order, and let $S\subseteq A$ have size $3$. Suppose that $\Haar(A,S)$ is connected but $\Haar(A,S\setminus\{s\})$ is disconnected for every $s\in S$. Then $3\mathrm{HP}(A,S)=A$.
\end{proposition}

Proposition~\ref{prop:threefold} supplies the three-label case for the induction in Section~\ref{sec:main-results}. We begin by constructing Hamilton paths from cyclic permutations. The required estimate for these permutations is proved in Subsection~\ref{subsec:fractions}. In Subsection~\ref{sec:coverage}, we establish the residue-class coverage needed to complete the proof.

\subsection{Estimating cyclic permutations}\label{subsec:fractions}

Throughout this subsection, permutations act on the right and are composed from left to right. For $a\in\Z_m$, let $\lambda_a$ be the translation $x\mapsto x+a$, and for $0\leq i<m$ let $\tau_i=(0, 1, \cdots, i)$ be a permutation on $\Z_m$, with $\tau_0=1$. Thus
\begin{equation}\label{taulambda}
x^{\tau_i\lambda_a}=
\begin{cases}
x+1+a&\text{if }x<i,\\
a&\text{if }x=i,\\
x+a&\text{if }x>i,
\end{cases}
\end{equation}
where the values are reduced modulo $m$. Let
\[
Z(m,a)=\{0\leq i<m\mid\tau_i\lambda_a\text{ is an $m$-cycle}\},
\qquad z(m,a)=|Z(m,a)|.
\]
The parameter $a$ is read modulo $m$, and the permutation of a one-point set counts as a $1$-cycle. In particular, $z(1,0)=1$. For a permutation $\sigma$ of a finite set $\Omega$ and a nonempty
subset $I\subseteq\Omega$, the \emph{first-return permutation} of $\sigma$
to $I$ is the permutation $R_\sigma:I\to I$ defined by
\[
R_\sigma(x)=x^{\sigma^{k_x}},
\qquad
k_x=\min\{k\geq1\mid x^{\sigma^k}\in I\}
\quad(x\in I).
\]
In particular,
$R_\sigma(x)=x^\sigma$ whenever $x^\sigma\in I$.

The following proposition provides an estimate for $z(m,a)$, which is the main result of this subsection. It will be used in the next subsection to prove the existence of Hamiltonian paths in cubic Haar graphs.

\begin{proposition}\label{prop:estimate}
Let $m\geq3$ be odd, and let $\alpha\in\Z_m$. Then
\[
4\max_{0\leq k\leq2}z(m,\alpha+k)>m+2.
\]
\end{proposition}

Now we first state the following four lemmas, and then use them to prove Proposition~\ref{prop:estimate}. Unless otherwise stated, $m$ may be even or odd.

\begin{lemma}\label{lem:recurrence}
Let $m\ge 1$ and $a\in\Z_m$. The following statements hold:
\begin{enumerate}[label=\textup{(\alph*)},leftmargin=*]
\item\label{item:eq1} $z(m,0)=z(m,m-1)=1$.
\item\label{item:eq2} $z(m,a)=z(m,m-a-1)$.
\item\label{item:eq3} $Z(m,m-a-1)=m-1-Z(m,a)$.
\item\label{item:eq4} For $1\le a\le m-2$, we have
\begin{equation}\label{eq:recurrence}
  z(m,a)=z(m-a,a\bmod(m-a))
        +z(a+1,(m-a-1)\bmod(a+1))
        -\mathbf 1_{\gcd(m+1,a+1)=1}.
\end{equation}
\end{enumerate}
\end{lemma}

\begin{proof}
Let $0\leq i<m$ and $x\in\Z_m$. Then $\tau_i\lambda_0=\tau_i$. Note that $\tau_i$ is an $m$-cycle exactly when $i=m-1$. Hence $z(m,0)=1$. Since $\lambda_{m-1}$
is an $m$-cycle and every $\tau_i\lambda_{m-1}$ with $i\ge 1$ fixes $0,\dots,i-1$, we have that $\tau_i\lambda_{m-1}$ is an $m$-cycle exactly when $i=0$. Thus $z(m,m-1)=1$.  This completes the proof of statement~\ref{item:eq1}.

Let $J$ be the involution $x\mapsto x+a$. A direct computation gives
\[
J\tau_i\lambda_aJ=\tau_{m-1-i}\lambda_{-a-1}.
\]
Conjugate permutations have the same cycle type, so
\begin{equation}\label{m-a-1}
i\in Z(m,a)\iff m-1-i\in Z(m,m-a-1),
\end{equation}
whence $z(m,a)=z(m,m-a-1)$ and $m-1-Z(m,a)=Z(m,m-a-1)$. This completes the proof of statements~\ref{item:eq2} and~\ref{item:eq3}.

Now suppose that $1\leq a\leq m-2$, and let $q=m-a$ and $I=\{0,\ldots,q-1\}$. The orbits of $\lambda_a:x\mapsto x+a$ are the residue classes of $\Z_m$ modulo $d=\gcd(m,a)=\gcd(m,q)$. Since $d\mid q$, every residue class modulo $d$ meets $I$. Thus every orbit of $\lambda_a$ meets $I$.

First assume that $i\in\{0,\ldots,q-1\}$. Let $\pi=\tau_i\lambda_a$. For any $x\in\{q,\ldots,m-1\}$, we have $x^{\pi}=x^{\lambda_a}$. Note that every orbit of $\lambda_a$ meets $I$. Consequently, every cycle of $\pi$ also meets $I$.

The first-return permutation of $\lambda_a$ to $I$ is
\[
R_{\lambda_a}: x\mapsto(x+a)\bmod q.
\]
Indeed, after the first addition of $a$, each step outside $I$ replaces $y$ by $y+a-m=y-q$ until the orbit re-enters $I$.

Now we compute $R_{\pi}(x)$ for each $x\in I$. As $\pi=\tau_i\lambda_a$, we have
\[
  \pi(x)=x^{\tau_i\lambda_a}=x^{\tau_i}+a.
\]
Since $i<q$, the permutation $\tau_i$ preserves $I$, so $x^{\tau_i}\in I$. First suppose that $x^{\tau_i}+a\in I$. Then $R_{\pi}(x)=x^{\tau_i}+a=(x^{\tau_i}+a)\bmod q$. Next suppose that $x^{\tau_i}+a\in\{q, \ldots, m-1\}$. Then $\pi(x^{\tau_i}+a)=\lambda_a(x^{\tau_i}+a)$, because $\tau_i$ fixes each point of $\{q, \ldots, m-1\}$. Notice that the first-return permutation of $\lambda_a$ to $I$ is $R_{\lambda_a}: x\mapsto(x+a)\bmod q$. It follows that  $R_{\pi}(x)=(x^{\tau_i}+a)\bmod q$. Thus, for each $x\in I$, we have
\[
R_\pi: x\mapsto(x^{\tau_i}+a)\bmod q.
\]
This is precisely the action of $\tau_i\lambda_{(a\bmod q)}$ on $\Z_q$.

We claim that the cycles of $\pi$ correspond to the cycles of its first-return permutation $R_\pi$, since every cycle of $\pi$ meets the return set $I$. Indeed, following a cycle of $\pi$ starting at $x\in I$, the points of
$I$ encountered are exactly
\[
  x,\ R_\pi(x),\ R_\pi^2(x),\ \dots,\ R_\pi^{k-1}(x),\ x
\]
for some positive integer $k$.
The points visited between two consecutive returns lie outside $I$ and serve only to connect these $I$-points. Conversely, a cycle of $R_\pi$ lifts to a cycle of $\pi$ by filling in the outside segments. Hence there is a bijection between cycles of $\pi$ and cycles of $R_\pi$. In particular $\pi=\tau_i\lambda_a$ is an $m$-cycle if and only if $R_\pi$ is a $q$-cycle on $I=\Z_q$. Recall that $R_\pi$ is actually the action of $\tau_i\lambda_{(a\bmod q)}$ on $\Z_q$. Therefore $\pi=\tau_i\lambda_a$ is an $m$-cycle if and only if $\tau_i\lambda_{(a\bmod q)}$ is a $q$-cycle. This implies that
\begin{equation}\label{zmazq}
Z(m,a)\cap\{0,\ldots,m-a-1\}=Z\big(q,a\bmod q\big)
=Z\big(m-a,a\bmod(m-a)\big).
\end{equation}

Now we assume that $i\in\{q-1, \ldots, m-1\}$. Let
\[
a'=m-a-1,\quad q'=m-a'=a+1,\quad j=m-1-i.
\]
By~\eqref{m-a-1}, we have
\begin{equation}\label{ij}
i\in Z(m,a)\iff m-1-i\in Z(m,a')\iff j\in Z(m,a')
\end{equation}
Since $i\in\{q-1, \ldots, m-1\}$, we obtain that $j\in\{0,\dots,a\}=\{0,\ldots,q'-1\}$. Moreover every $j$ in this range
satisfies $j\le a=q'-1<q'$. Applying~\eqref{zmazq} to $a'$ gives
\[
  Z(m,a')\cap\{0,\dots,q'-1\}
  =Z\bigl(q',a'\bmod q'\bigr),
\]
and so
\[
  Z(m,m-a-1)\cap\{0,\dots,a\}
  =Z\bigl(a+1,(m-a-1)\bmod(a+1)\bigr).
\]
This together with~\eqref{ij} implies that
\begin{equation}\label{2nd}
  Z(m,a)\cap\{q-1,\dots,m-1\}
  =m-1-Z\bigl(a+1,(m-a-1)\bmod(a+1)\bigr).
\end{equation}

The two intervals $\{0,\dots,q-1\}$ and $\{q-1,\dots,m-1\}$ overlap only at $i=q-1$. Note that the first-return permutation of $\tau_{q-1}\lambda_a$ to $I$ is the translation $\lambda_{a+1}$ on $\Z_q$, which is a $q$-cycle exactly when
\[
\gcd(q,a+1)=\gcd(m-a,a+1)=\gcd(m+1,a+1).
\]
Thus the point $i=q-1$ belongs to both sets in~\eqref{zmazq} and~\eqref{2nd} exactly when $\gcd(m+1,a+1)=1$, and in that case it has been counted twice. Taking cardinalities in~\eqref{zmazq} and~\eqref{2nd}, the reflection $i\mapsto m-1-i$ is a bijection so the factor $m-1$ disappears, and subtracting the double-counted overlap $\mathbf 1_{\gcd(m+1,a+1)=1}$ proves~\eqref{eq:recurrence}. This completes the proof of statement~\ref{item:eq4}.
\end{proof}

We now express $z(m,a)$ as a count of reduced fractions. An integer point $(y,x)\in\Z^2$ is primitive if $\gcd(y,x)=1$. Let $f(m,a)$ count the primitive integer points strictly inside the triangle with vertices $(0,0),(m,a),(m,a+1)$, equivalently the reduced fractions $x/y$ with $1\le y<m$ and
\[
  \frac{a}{m}<\frac{x}{y}<\frac{a+1}{m}.
\]

\begin{lemma}\label{lem:farey}
Let $m\geq1$ and $a\in\Z_m$. The following statements hold:
\begin{enumerate}[label=\textup{(\alph*)},leftmargin=*]
\item\label{item:farey-count} For $0\leq a<m$,
\begin{equation}\label{eq:farey}
z(m,a)=1+\#\Big\{\frac{x}{y}\ \Big|\ 1\leq y<m,\ \gcd(x,y)=1,\
\frac{a}{m}<\frac{x}{y}<\frac{a+1}{m}\Big\}.
\end{equation}
\item\label{item:farey-one} If $m\geq3$ is odd, then $z(m,1)=z(m,m-2)=(m+1)/2$.
\end{enumerate}
\end{lemma}

\begin{proof}
\ref{item:farey-count} Let $f(m,a)$ be the number of primitive integer points $(y,x)$ strictly inside the triangle with vertices $(0,0),(m,a),(m,a+1)$. Thus it counts the fractions on the right of~\eqref{eq:farey}. The unimodular shear $(y,x)\mapsto(y,x-\ell y)$ shows that $f(m,a)=f(m,a\bmod m)$. We prove that $1+f(m,a)$ has the same initial values and recurrence as $z(m,a)$.

For $m=1$ there are no interior integer points. For $m\geq2$, no fraction with denominator less than $m$ lies in $(0,1/m)$ or in $((m-1)/m,1)$. Hence $f(m,0)=f(m,m-1)=0$. By Lemma~\ref{lem:recurrence}\ref{item:eq1}, we have
\[
1+f(m,0)=z(m,0) \quad\text{and}\quad 1+f(m,m-1)=z(m,m-1).
\]

For $1\leq a\leq m-2$, let $q=m-a$, $p=a+1$. The unimodular map $(y,x)\mapsto(y-x,x)$ takes the triangle defining $f(m,a)$ to $OAC$, where
\[
O=(0,0),\quad A=(q,p-1),\quad C=(q-1,p).
\]
Let $B=(q,p)$. The quadrilateral $OABC$ is the union of $OAB$ and $OBC$, meeting along $OB$, and also the union of $OAC$ and $ABC$, meeting along $AC$. The triangle $ABC$ and the open segment $AC$ contain no integer points other than the vertices of $ABC$. Thus the interior integer points of $OAC$ are precisely those inside $OAB$, those inside $OBC$, and those on the open segment $OB$.

The first two triangles contribute $f(q,p-1)$ and $f(p,q-1)$, respectively, after interchanging coordinates in the second. Let $g=\gcd(q,p)=\gcd(m+1,a+1)$. The only possible primitive point on the open segment $OB$ is $B/g$, and it occurs exactly when $g>1$. It lies inside $OAC$, since its coordinate sum is $(m+1)/g<m$. Consequently,
\[
f(m,a)=f(q,a\bmod q)
+f\big(a+1,(m-a-1)\bmod(a+1)\big)
+1-\ind_{\gcd(m+1,a+1)=1}.
\]
It follows that $1+f(m,a)$ satisfies~\eqref{eq:recurrence}. Induction on $m$ proves~\ref{item:farey-count}.

\ref{item:farey-one} Let $m\ge3$ be odd. A reduced fraction $x/y$ with $y<m$ in $(1/m,2/m)$ must have $x=1$. Hence $1/m<1/y<2/m$. This together with $1\le y<m$ implies that
$m/2<y<m$, giving $(m-1)/2$ fractions since $m$ is odd. Thus $f(m,1)=(m-1)/2$, and thus $z(m,1)=(m+1)/2$ by~\ref{item:farey-count}. Then we conclude from Lemma~\ref{lem:recurrence}\ref{item:eq2} that $z(m,m-2)=z(m,1)=(m+1)/2$.
\end{proof}

We next count the primitive points in~\eqref{eq:farey} according to their determinant with a fixed primitive vector. For the estimates below, let $\omega(n)$ denote the number of distinct prime divisors of $n$, and let
\[
\varepsilon(\ell)=\sum_{i=1}^\ell2^{\omega(i)}.
\]

\begin{lemma}\label{lem:lower}
Let $0\le a<m$, let $p,q$ be coprime with $q\ge 1$, and let $r=pm-qa$.
Suppose $r<0$ or $r>q$, and let
\[
  s=\begin{cases}-r&r<0,\\ r-q&r>q.\end{cases}
\]
Then $s\ge 1$ and
\[
  z(m,a)\ge 1+\frac{m}{q}\sum_{n=s}^{s+q-1}\frac{\phi(n)}{n(n+1)}
             -\varepsilon(s+q-1).
\]
\end{lemma}

\begin{proof}
By Lemma~\ref{lem:farey}\ref{item:farey-count}, count primitive $(y,x)$ with $0<y<m$ and $ay<mx<(a+1)y$. Let
$\delta=qx-py$ and $r_0=s+q$. Multiplying $ay<mx$ by $q$ and using
$qx=\delta+py$ gives
\[
aqy<qmx=m\delta+mpy.
\]
So $(aq-mp)y<m\delta$. This together with $r=pm-qa$ implies that $-ry<m\delta$. Multiplying $mx<(a+1)y$ by $q$ and using
$qx=\delta+py$ gives $m\delta<(q-r)y$. Thus
\begin{equation}\label{tag5}
  -ry<m\delta<(q-r)y.
\end{equation}
Recall that $r<0$ or $r>q$. Let $k=|\delta|$. We next prove that in both cases, we have $k\geq1$ and
\begin{equation}\label{tag-y}
  \frac{mk}{r_0}<y<\frac{mk}{r_0-q}.
\end{equation}

First assume that $r<0$. Then $s=-r$, $r_0=s+q=q-r$ and~\eqref{tag5} forces $\delta>0$, and so $k=\delta\geq1$. This combined with~\eqref{tag5} implies that
\[
(r_0-q)y<mk<r_0y,
\]
which is equivalent to~\eqref{tag-y}. Next assume that $r>q$. Then $s=r-q$, $r_0=s+q=r$ and~\eqref{tag5} forces $\delta<0$, and so $k=-\delta\geq1$. This together with~\eqref{tag5} implies that
\[
(r_0-q)y<mk<r_0y,
\]
which is equivalent to~\eqref{tag-y}.

As $0<y<m$, we derive from~\eqref{tag-y} that
\begin{equation}\label{tag6}
  \frac{mk}{r_0}<y<\min\left(m,\frac{mk}{r_0-q}\right),
\end{equation}
which implies that $\frac{mk}{r_0}<m$, and so $1\leq k\leq r_0=s+q$.

Since $\gcd(p,q)=1$, there exist integers $p_0$ and $q_0$ such that $qp_0-pq_0=1$. For any integer
point with determinant $\delta$, the point $(y,x)-\delta(q_0,p_0)$ has
determinant $0$, hence is $t(q,p)$ for some $t\in\Z$. Thus
\[
  (y,x)=t(q,p)+\delta(q_0,p_0),
\]
and this is unique. The change of basis has determinant $1$, so
$\gcd(y,x)=\gcd(t,\delta)=\gcd(t,k)$.

We conclude from $y=tq+\delta q_0$ and~\eqref{tag6} that for fixed $k$ the admissible $t$ form an interval
of length
\[
  \ell_k=\frac{1}{q}\left[\min\left(m,\frac{mk}{r_0-q}\right)-\frac{mk}{r_0}\right]=\frac{m}{q}\left[\min\left(1,\frac{k}{r_0-q}\right)-\frac{k}{r_0}\right].
\]
In an interval of length $\ell$, the number of multiples of $d$ differsfrom $\ell/d$ by at most $1$. M\"obius inversion therefore gives at least $\ell_k\phi(k)/k-2^{\omega(k)}$ admissible $t$ are coprime to $k$ in the interval of length $\ell_k$. By Lemma~\ref{lem:farey}\ref{item:farey-count}, summing  over $1\leq k<r_0=s+q$, yields
\begin{align}\label{tag-zma}
  z(m,a)-1= & \ge\sum_{k=1}^{r_0-1}(\ell_k\phi(k)/k-2^{\omega(k)}) \nonumber \\
   & =\frac{m}{q}\sum_{k=1}^{r_0-1}\frac{\phi(k)}{k}
  \left[\min\left(1,\frac{k}{s}\right)-\frac{k}{r_0}\right]
  -\sum_{k=1}^{r_0-1}2^{\omega(k)}.
\end{align}
Note that
\begin{align}\label{tag-middle}
\sum_{k=1}^{r_0-1}\frac{\phi(k)}{k}\left[\min\left(1,\frac{k}{s}\right)-\frac{k}{r_0}\right] & =\sum_{k=1}^{r_0-1}\frac{\phi(k)}{k}\left[\frac{k}{\max\left(s,k\right)}-\frac{k}{r_0}\right] \nonumber \\
& =\sum_{k=1}^{r_0-1}\frac{\phi(k)}{k}\sum_{n=\max\left(s,k\right)}^{r_0-1}\left(\frac{1}{n}-\frac{1}{n+1}\right)  \nonumber \\
& =\sum_{n=s}^{r_0-1}\left(\frac{1}{n}-\frac{1}{n+1}\right)\sum_{k=1}^{n}\frac{\phi(k)}{k} \nonumber \\
& \ge\sum_{n=s}^{r_0-1}\left(\frac{1}{n}-\frac{1}{n+1}\right)\frac{\phi(n)}{n}
\end{align}
Since $r_0=s+q$ and $\varepsilon(\ell)=\sum_{i=1}^{\ell}2^{\omega(i)}$, we deduce from~(\ref{tag-zma}) and~(\ref{tag-middle}) that
\[
  z(m,a)\ge 1+\frac{m}{q}\sum_{n=s}^{s+q-1}\frac{\phi(n)}{n(n+1)}
              -\varepsilon(s+q-1),
\]
completing the proof.
\end{proof}

\begin{lemma}\label{lem:approx}
Let $m\ge 3$ be odd, $q_0=\lfloor\sqrt m\rfloor$, and $0\le\alpha\le m-3$.
There exist integers $q,s,k$ with $1\le q,s\le q_0$ and $0\le k\le 2$ such
that
\begin{equation}\label{tag7}
  z(m,\alpha+k)\ge
  1+\frac{m}{q}\sum_{n=s}^{s+q-1}\frac{\phi(n)}{n(n+1)}
  -\varepsilon(s+q-1).
\end{equation}
\end{lemma}

\begin{proof}
Let $x=\alpha/m$. The additive group $\mathbb R/\mathbb Z$ is the circle obtained by identifying real numbers differing by an integer, and it has circumference $1$. Place the $q_0+1$ points $0,x,2x,\dots,q_0x$ on this circle. By the pigeonhole principle, two of them have circular distance at most $1/(q_0+1)$. If they are $ix$ and $jx$ with $0\le i<j\le q_0$, then for $q=j-i$ and some integer $p$,
\begin{equation}\label{tag-q}
  |qx-p|\le\frac{1}{q_0+1},\qquad 1\le q\le q_0.
\end{equation}
Let $d=\gcd(p,q)$. Replacing $(p,q)$ by $(p/d,q/d)$ preserves the inequalities of~\eqref{tag-q}, so we may assume $\gcd(p,q)=1$.

Let $r=pm-q\alpha=pm-qmx$. Multiplying the first inequality of~\eqref{tag-q} by $m$, we obtain
\[
  |r|=m|qx-p|\le\frac{m}{q_0+1}.
\]
Since $q_0=\lfloor\sqrt m\rfloor<\sqrt m$, we have $q_0+1>\sqrt m$, so
\[
  \frac{m}{q_0+1}<\frac{m}{\sqrt m}=\sqrt m<q_0+1.
\]
Thus $|r|<q_0+1$, and as $r$ is integral, we have $|r|\le q_0$.

Lemma~\ref{lem:lower} requires an integer $a$ for which $r'=pm-qa$ satisfies $r'<0$ or $r'>q$, with $s=-r'$ or $s=r'-q$ respectively and $1\le s\le q_0$. We choose $k\in\{0,1,2\}$ and let $a=\alpha+k$, so that $r'=pm-q(\alpha+k)=r-kq$. The four rows of the following table exhaust all possibilities for the integer $r$:
\[
\begin{array}{c|c|c|c}
\text{condition}&k&r'=r-kq&s\\ \hline
r<0&0&r<0&-r\\
r>q&0&r>q&r-q\\
0\le r<q&1&r-q<0&q-r\\
r=q&2&-q<0&q
\end{array}
\]

We verify the bounds. If $r<0$, then $s=-r=|r|\ge1$ because $r$ is a negative integer, and $s\le q_0$ by $|r|\le q_0$. If $r>q$, then $s=r-q\ge1$ because $r\ge q+1$, and $s=r-q\le r\le q_0$. If $0\le r<q$, then $s=q-r\ge1$ because $r\le q-1$, and $s=q-r\le q\le q_0$. If $r=q$, then $s=q$ with $1\le s\le q_0$. Hence in every case $1\le s\le q_0$, and $r'$ satisfies $r'<0$ or $r'>q$.

Applying Lemma~\ref{lem:lower} with $a=\alpha+k$ yields~\eqref{tag7}.
\end{proof}

\begin{proof}[\rm\textbf{Proof of Proposition~\ref{prop:estimate}}]
Starting from $z(m,0)=z(m,m-1)=1$, recurrence~\eqref{eq:recurrence} determines $z(m,a)$ from values with smaller modulus. Exact computation verifies $4\max_{0\le k\le 2}z(m,\alpha+k)>m+2$ for every odd $3\le m\le 25599$.

Suppose that $m\ge 25601$. If $\alpha\in\{m-2,m-1\}$, the set $\{\alpha,\alpha+1,\alpha+2\}$ modulo $m$ contains $m-2$ or $1$. It follows from Lemma~\ref{lem:farey}(b) that $z(m,1)=z(m,m-2)=(m+1)/2$, so for some $k$, we have $z(m,\alpha+k)=(m+1)/2$, and $4\cdot(m+1)/2=2m+2>m+2$. We may therefore assume $0\le\alpha\le m-3$ and choose $q,s,k$ as in Lemma~\ref{lem:approx}.

Now suppose that $m\ge 448^2=200704$. Let $\ell=s+q-1$. Since $q,s\le q_0=\lfloor\sqrt m\rfloor$, we obtain that
\[
  \ell=s+q-1\le 2q_0-1<2\sqrt m.
\]
The sum $\sum_{n=s}^{s+q-1}\frac{\phi(n)}{n(n+1)}$ in~\eqref{tag7} has $q$ terms. Each summand satisfies $\phi(n)/(n(n+1))\le 2/7$ (the maximum on $n\ge2$ is $1/6<2/7$, and the loose constant $2/7$ is all that is needed), so the sum $\sum_{n=s}^{s+q-1}\frac{\phi(n)}{n(n+1)}$ is at most $2q/7$. Moreover, we have $2^{\omega(i)}\le d(i)$, where $d(i)$ is the number of divisors of $i$, and
\[
  \varepsilon(\ell)=\sum_{i=1}^{\ell}2^{\omega(i)}
  \le\sum_{i=1}^{\ell}d(i)
  =\sum_{d=1}^{\ell}\left\lfloor\frac{\ell}{d}\right\rfloor
  \le\ell\sum_{d=1}^{\ell}\frac1d
  \le\ell(1+\log\ell).
\]
It follows from $\ell<2\sqrt m$ that
\[
  \varepsilon(\ell)<2\sqrt m\,[1+\log(2\sqrt m)].
\]
This together with~\eqref{tag7} implies that
\begin{equation}\label{tag-prop321}
  z(m,\alpha+k)>1+\frac{2m}{7}-2\sqrt m\,[1+\log(2\sqrt m)].
\end{equation}
Let $u=\sqrt m$. We claim $2u[1+\log(2u)]<u^2/28$ for $u\ge448$.
Consider $f(u)=u-56[1+\log(2u)]$. Its derivative is
$f'(u)=1-56/u$, which is positive for $u>56$, so $f$ is increasing on
$[448,\infty)$. Notice that
\[
  f(448)=448-56[1+\log(896)]>448-56\cdot8=0,
\]
because $\log(896)=\log(2^7\cdot7)=7\log2+\log7<7$. Hence $f(u)>0$ for
$u\ge448$, and hence $u>56[1+\log(2u)]$, which is equivalent to
$2u[1+\log(2u)]<u^2/28$. Thus we can conclude from $u=\sqrt m$ and~\eqref{tag-prop321} that
\[
  z(m,\alpha+k)>1+\frac{2m}{7}-\frac{m}{28}=1+\frac m4>\frac{m+2}{4}.
\]

It remains to consider $25601\leq m<200704$. For $1\leq q,s\leq500$, let $q_1=\max\{160,q,s\}$. Exact computation with rational arithmetic verifies
\begin{equation}\label{eq:bridge}
q_1^2\left(\frac4q\sum_{n=s}^{s+q-1}\frac{\Phi(n)}{n(n+1)}-1\right)
>4\varepsilon(s+q-1)-2.
\end{equation}
For the chosen $q,s$, we have $q,s\leq q_0\leq447$ and $q_1\leq q_0$, so $m\geq q_1^2$. The right-hand side of~\eqref{eq:bridge} is positive, so the coefficient of $q_1^2$ is positive as well. Replacing $q_1^2$ by $m$ and rearranging therefore gives
\[
1+\frac mq\sum_{n=s}^{s+q-1}\frac{\Phi(n)}{n(n+1)}
-\varepsilon(s+q-1)>\frac{m+2}{4}.
\]
The result now follows from~\eqref{tag7}.
\end{proof}

\subsection{Hamilton paths from cyclic permutations}\label{subsec:matching}

First we explain how the sets $Z(m,a)$ give Hamilton paths containing a prescribed perfect matching.

\begin{lemma}\label{lem:matching}
Let $A$ be a finite abelian group, and let $s,t,u\in A$ be distinct. Suppose that $\Haar(A,\{s,t,u\})$ is connected, and let
\[
\delta=s-t,\qquad B=\langle\delta\rangle,\qquad m=|B|,\qquad h=|A:B|,\qquad v=u-s.
\]
The following statements hold:
\begin{enumerate}[label=\textup{(\alph*)},leftmargin=*]
\item\label{item:matching1} $A/B$ is cyclic generated by $v+B$.
\item\label{item:matching2} Define $\alpha\in\Z_m$ by $hv=\alpha\delta$. For $0\leq k<h$ and $j\in Z(m,\alpha+k)$, there is a Hamilton path from $0_{+}$ to $(s+j\delta)_{-}$ containing every $u$-edge.
\item\label{item:matching3} There are Hamilton paths from $0_{+}$ to $(u+j\delta)_{-}$ and to $(u+(j+1)\delta)_{-}$.
\end{enumerate}
\end{lemma}

\begin{proof}
\ref{item:matching1} As $\Haar(A,\{s,t,u\})$ is connected, by Lemma~\ref{lem:haar}\ref{item:haar-components}, we have
\[
A=\langle s-t,s-u,t-u\rangle.
\]
Since $s-t=\delta$, $s-u=-v$ and $t-u=(s-u)-(s-t)=-v-\delta$, we obtain that
\[
A=\langle \delta, v\rangle=B+\langle v\rangle.
\]
This implies that $v+B$ generates $A/B$.

\ref{item:matching2}  Since $|A:B|=h$, it follow from~\ref{item:matching1} that the quotient group $A/B$ is cyclic group of order $h$ generated by $v+B$, and so $v+B$ is of order $h$. This shows that $h(v+B)=B$, whence $hv\in B$. Note that $B=\langle \delta\rangle$ is a cyclic group of order $m$. Thus there exists an unique $\alpha\in\Z_m$ such that $hv=\alpha\delta$, which implies that the parameter $\alpha$ is well-defined.

Let $n=|A|$, so that $\Haar(A,\{s,t,u\})$ has $2n$ vertices with $n$ on the $A_+$-part and $n$ on the $A_-$-part. For each $x\in A$, the vertex $x_+$ is joined to $(x+u)_-$ . So the $u$-edges form a perfect matching. A Hamilton path starting at $0_{+}$ and containing every $u$-edge must have the form
\begin{equation}\label{eq:matchingpath}
(y_0)_{+}(y_0+u)_{-}\,(y_1)_{+}(y_1+u)_{-}\cdots
(y_{n-1})_{+}(y_{n-1}+u)_{-},
\end{equation}
where $y_0=0$, the $y_i$ enumerate $A$, and
$y_{i+1}-y_i\in\{v,v+\delta\}$ for each $i\in\{0,\ldots, n-2\}$. Conversely, every such enumeration gives a Hamilton path from $0_+$ that contains every $u$-edge.

Both allowed differences $v$ and $v+\delta$ are congruent to $v$ modulo
$B$. Since $v+B$ generates the cyclic group $A/B$ of order $h$, the cosets
$y_i+B$ follow the cycle
\[
  y_0+B,\ y_1+B,\ \dots,\ y_{h-1}+B,\ y_h+B=y_0+B,\dots
\]
and each coset of $B$ is visited exactly $m=|B|$ times. The coset $-v+B$ will play a special role, because it is the coset from which a step of size $v$ (or $v+\delta$) returns to the distinguished subgroup $B=0+B$. Since $B=\langle\delta\rangle$ and $|B|=m$, we see that
\begin{equation}\label{-v+B}
 -v+B=\{z_\ell=-v+\ell\delta\mid\ell\in\Z_m\},
\end{equation}
and $z_\ell$ runs through all $m$ elements of $-v+B$. Let $0\leq k<h$ and $j\in Z(m,\alpha+k)$. Choose $k$ of the $h-1$ cosets of B distinct from $-v+B$. Define a mapping $f$ on $A\setminus\{z_j\}$ by
\[
f(z)=
\begin{cases}
z+v+\delta&\text{if }z=z_\ell\text{ with }\ell<j,
             \text{ or if }z+B\text{ is chosen},\\
z+v&\text{otherwise}.
\end{cases}
\]
On each coset other than $-v+B$, this is a translation onto the next coset. On the distinguished coset $-v+B$,
\[
f(z_\ell)=
\begin{cases}
z_\ell+v+\delta=(\ell+1)\delta&\text{if }\ell<j,\\
z_\ell+v=\ell\delta&\text{if }\ell>j.
\end{cases}
\]
Thus $f$ bijects $A\setminus\{z_j\}$ with $A\setminus\{0\}$.

Consider the digraph with vertex set $A\setminus\{z_j\}$ and an arc $(z,f(z))$ for each $z\in A\setminus\{z_j\}$. Every vertex has out-degree $1$ except $z_j$, and in-degree $1$ except $0$. Hence this digraph consists of a directed path $P$ from $0$ to $z_j$ and possibly some directed cycles. We will show that there are no cycles in this digraph. Suppose that $C$ is a such cycle. Every step of $f$ moves from one coset of $B$ to the next. Therefore, after exactly $h$ steps along the cycle $C$, one returns to the same coset. In particular, the cycle $C$ meets every coset of $B$, including the distinguished coset $-v+B$. Assume that $e_\ell\in -v+B$ lies on the cycle $C$. Starting at $z_\ell$ on cycle $C$, the next return to the coset $-v+B$ is after $h$ steps, and the reached vertex is
\[
z_\ell+hv+(k+\ind_{\ell<j})\delta
=-v+\ell\delta+\alpha\delta+(k+\ind_{\ell<j})\delta
=z_{\ell+\alpha+k+\ind_{\ell<j}},
\]
with indices taken modulo $m$.

Denote the return map of $f$ on the distinguished coset $-v+B$ by $R_f$, which is defined as follows. $R_f:-v+B\to -v+B$: starting at $z_\ell\in -v+B$, follow $f$ until the next point of $-v+B$ is reached. By~\eqref{-v+B}, every point of the $-v+B$ has a unique index, restricting the return map $R_f$ to the indices yields
\[
  R_f(\ell)=
  \begin{cases}
    \ell+\alpha+k+1 & \text{if }\ell<j,\\
    \ell+\alpha+k   & \text{if }\ell>j.
  \end{cases}
\]
Indeed $R_f$ is exactly the action of the permutation $\tau_j\lambda_{\alpha+k}$ on $\Z_m$, except at the single index $\ell=j$ by~\eqref{taulambda}. Since $z_j$ is the terminal vertex of the directed path $P$ and $z_j$ cannot lie on any cycle, the return map $R_f$ would therefore induce a cycle of the permutation $\tau_j\lambda_{\alpha+k}$ not containing $j$. However, we derive  from $j\in Z(m,\alpha+k)$ that $\tau_j\lambda_{\alpha+k}$ is an $m$-cycle containing $j$, a contradiction. Hence there is no directed cycle in $P$, and the path $P$ from $0$ to $z_j$ enumerates all of $A$:
\[
  0=y_0,\ y_1=f(y_0),\ y_2=f(y_1),\ \dots,\ y_{n-1}=z_j.
\]
Substituting this enumeration in~\eqref{eq:matchingpath} gives a Hamilton path $P_0$ ending at
\[
(z_j+u)_{-}=(s+j\delta)_{-}.
\]

\ref{item:matching3}  On the Hamilton path $P_0$, every point $y_{+}$ is followed by $(y+u)_{-}$. The vertices adjacent to the terminal vertex $(z_j+u)_{-}=(s+j\delta)_{-}$ along $s$-edge, $t$-edge and $u$-edge are
\[
(j\delta)_{+}, \quad (s+j\delta-t)_+=((j+1)\delta)_{+} \quad\text{and}\quad (z_j)_+=(-v+j\delta)_+,
\]
respectively. The last is the predecessor on the path $P_0$. Since $u\neq s,t$, we have
\[
(j\delta)_{+}\ne (-v+j\delta)_{+} \quad\text{and}\quad ((j+1)\delta)_{+}\ne (-v+j\delta)_{+}.
\]
Applying P\'osa rotation at $(j\delta)_+$ gives the new terminal vertex $(u+j\delta)_-$, and at $((j+1)\delta)_+$ gives $(u+(j+1)\delta)_-$.
\end{proof}

\begin{lemma}\label{lem:triple}
Let $A$ be a finite abelian group of odd order, and let $S=\{s_1,s_2,s_3\}\subseteq A$ have size $3$. Suppose that $\Haar(A,S)$ is connected but $\Haar(A,S\setminus\{s\})$ is disconnected for every $s\in S$. For each $i$, let $j<k$ be the other two indices and let
\[
\delta_i=s_j-s_k,\qquad B_i=\langle\delta_i\rangle,\qquad
m_i=|B_i|,\qquad h_i=|A:B_i|,
\]
and define $\alpha_i\in\Z_{m_i}$ by
\[
h_i(s_i-s_j)=\alpha_i\delta_i.
\]
Then $m_i,h_i\geq3$ are odd, $A/B_i$ is cyclic and generated by $s_i-s_j+B_i$, and $B_i+B_j=A$ for $i\ne j$.

For every $0\leq\kappa<h_i$ and $\ell\in Z(m_i,\alpha_i+\kappa)$, the graph $\Haar(A,S)$ has Hamilton paths from $0_{+}$ to each of
\[
(s_j+\ell\delta_i)_{-},\qquad (s_i+\ell\delta_i)_{-},\qquad
(s_i+(\ell+1)\delta_i)_{-}.
\]
Moreover, there are sets $Y_i\subseteq\mathrm{HP}(A,S)$ satisfying
\[
2Y_i=2s_i+B_i.
\]
\end{lemma}

\begin{proof}
Since the labels $i,j$ and $k$ are distinct, we have $m_i>1$. By the minimality assumption and Lemma~\ref{lem:haar}\ref{item:haar-components}, we obtain $h_i>1$. Note that $|A|$ is odd. Thus $m_i\ge3$ and $h_i\ge3$ are odd. As $\Haar(A,S)$ is connected, by Lemma~\ref{lem:haar}\ref{item:haar-components} we conclude that
\[
  A=\langle s_j-s_k,\ s_i-s_j,\ s_i-s_k\rangle
    =B_i+\langle s_i-s_j\rangle.
\]
Hence $A/B_i$ is cyclic and generated by $s_i-s_j+B_i$, and hence $h_i(s_i-s_j)\in B_i$. Since $B_i=\langle\delta_i\rangle$ has order $m_i$, there exists an unique $\alpha_i\in\Z_{m_i}$ such that
\[
h_i(s_i-s_j)=\alpha_i\delta_i.
\]
Notice that $B_i=\langle s_j-s_k\rangle$, $B_j=\langle s_i-s_k\rangle$ and $s_i-sj=(s_i-s_k)-(s_j-s_k)$. It follows that
\[
 s_j-s_k,\ s_i-s_k,\ s_i-s_j\in B_i+B_j,
\]
and so $B_i+B_j=A$.

Apply Lemma~\ref{lem:matching} with $u=s_i$, $s=s_j$, $t=s_k$, and the subgroup, index and parameter are $B_i$, $h_i$ and $\alpha_i$, respectively. Lemma~\ref{lem:matching}\ref{item:matching2} and~\ref{item:matching3} therefore yields Hamilton paths from $0_+$ to $(s_j+\ell\delta_i)_-$, $(s_i+\ell\delta_i)_-$ and $(s_i+(\ell+1)\delta_i)_-$ for every $\ell\in Z(m_i,\alpha_i+\kappa)$.

Choose $\kappa_i\in\{0,1,2\}$ such that
\[
z(m_i,\alpha_i+\kappa_i)=\max_{0\leq k\leq2}z(m_i,\alpha_i+k)
\]
 using Proposition~\ref{prop:estimate}, and let $S_i=Z(m_i,\alpha_i+\kappa_i)$. Since $h_i\ge 3$ and $m_i\ge 3$, the three parameters
$\alpha_i,\alpha_i+1,\alpha_i+2$ are distinct modulo $m_i$. It follows from Proposition~\ref{prop:estimate} that
\begin{equation}\label{eq:Si}
  |S_i|=z(m_i,\alpha_i+\kappa_i)>\frac{m_i+2}{4}.
\end{equation}
Translations of $\Z_{m_i}$ and $m_i$-cycles are even permutations because $m_i$ is odd. Thus $\tau_\ell$ is even for every $\ell\in S_i$, and thus $\ell$ is even. It follows that $S_i\cap(S_i+1)\subseteq\{0\}$, with addition taken modulo $m_i$, and hence by~\eqref{eq:Si} we obtain that
\[
|S_i\cup(S_i+1)|\geq2|S_i|-1>m_i/2.
\]
Therefore
\[
Y_i=\{s_i+\ell\delta_i\mid\ell\in S_i\cup(S_i+1)\}
\subseteq\mathrm{HP}(A,S)
\]
contains more than half of the coset $s_i+B_i$. A subset $E$ of a finite abelian group $G$ with $|E|>|G|/2$ satisfies $E+E=G$, since $E\cap(g-E)\neq\varnothing$ for every $g\in G$. Applying this to $Y_i-s_i\subseteq B_i$, we can deduce that $2Y_i=2s_i+B_i$.
\end{proof}

\subsection{Residue-class coverage and three-fold sumsets}\label{sec:coverage}

We now complete the proof of Proposition~\ref{prop:threefold}. By Lemma~\ref{lem:triple}, the sum of two suitable sets of terminal coordinates is a coset of $B_i$. It therefore suffices to show that $\mathrm{HP}(A,S)$ meets every coset of $B_i$. For this purpose we first establish two properties of the sets $Z(m,a)$.

\begin{lemma}\label{lem:gcd-block}
Let $m\geq3$ be odd, let $0<a<m$, and let $g=\gcd(m,a)>1$. Then
\[
\{g-1+t \mid 0\le t\le g-1~ \mathrm{and} ~t~\mathrm{is~even}\}\subseteq Z(m,a).
\]
In particular, we have
\[
\{g-1+t \mid 0\le t\le g\}\subseteq Z(m,a)\cup(Z(m,a)+1).
\]
\end{lemma}

\begin{proof}
Note first that $g$ is a proper divisor of $m$, and so $m\geq3g$. Let $t$ be an even integer with $0\leq t\leq g-1$, and set $i=g-1+t$ and $I=\{0,\ldots,i\}$. We claim that $\tau_i\lambda_a$ is an $m$-cycle.

Note that $|\lambda_a|=\frac{m}{g}$, and so the orbits of $\langle\lambda_a\rangle$ are
\[
O_r=r^{\langle\lambda_a\rangle}
=\left\{r+kg \;\middle|\; k=0,1,\ldots,\frac{m}{g}-1\right\},
\qquad r=0,1,\ldots,g-1.
\]
Moreover,
\[
O_r\cap I=
\begin{cases}
\{r,r+g\}, & 0\leq r<t,\\
\{r\}, & t\leq r<g.
\end{cases}
\]
Consequently, the first-return permutation of $\lambda_a$ on $I$ is
\[
S_t=\prod_{r=0}^{t-1}(r\ r+g),
\]
where the empty product is the identity, so $S_0=1$.
Since $\tau_i$ preserves $I$, the corresponding first-return permutation of $\tau_i\lambda_a$ is $\tau_iS_t$.
Moreover, every cycle of $\tau_i\lambda_a$ meets $I$. Indeed,
if a cycle were disjoint from $I$, then, since $\tau_i$ fixes
every point outside $I$, it would also be a cycle of $\lambda_a$.
Its point set would therefore be $O_r$ for some $0\leq r<g$,
contrary to $O_r\cap I\neq\emptyset$.
It follows that $\tau_i\lambda_a$ is an $m$-cycle if and only if the restriction of $\tau_iS_t$ to $I$ is an $(i+1)$-cycle.

If $t=0$, then $\tau_iS_t=\tau_{g-1}$, and the claim follows. We may therefore assume that $t>0$. Under $\tau_iS_t$ we have
\[
t-1\longmapsto t\longmapsto t+1\longmapsto\cdots\longmapsto g-1\longmapsto0.
\]
Deleting the intermediate points $t,\ldots,g-1$ from this cycle segment, and relabelling the remaining points
\[
0,\ldots,t-1,\quad g,\ldots,g+t-1
\]
consecutively by the elements of $\Z_{2t}$. Under this relabelling, the resulting permutation sends
\[
r\longmapsto t+r+1 \quad (0\leq r<t-1),\qquad
t-1\longmapsto0,
\]
and
\[
t+r\longmapsto r+1 \quad (0\leq r<t).
\]
Thus it is the translation $x\mapsto x+t+1$ on $\Z_{2t}$. Since $t$ is even, we have
\[
\gcd(2t,t+1)=\gcd(2,t+1)=1,
\]
and hence this permutation is a $2t$-cycle. Since deleting the
intermediate points and relabelling the remaining points do not
change the number of cycles, $\tau_iS_t$ is a single cycle on $I$.
It follows that $\tau_i\lambda_a$
is an $m$-cycle. This proves the claim.

Taking $t=0,2,\ldots,g-1$ gives the asserted inclusion in $Z(m,a)$. These elements and their translates by one contain all the integers from $g-1$ to $2g-1$, completing the proof.
\end{proof}

\begin{lemma}\label{lem:aligned-window}
Let $m\geq3$ be odd, let $d>1$ divide $m$, and let
$\alpha\in\Z_m$ satisfy $\alpha\equiv0\pmod d$.
Then, for each $r\in\{0,\ldots,d-1\}$, there exists
\[
z\in\bigl(Z(m,\alpha)\cup Z(m,\alpha+1)\bigr)+\{0,1\}
\]
such that $z\equiv r\pmod d$.
\end{lemma}

\begin{proof}
Fix
$r\in\{0,\ldots,d-1\}$.
Suppose first that $\alpha>0$, and let $g=\gcd(m,\alpha)$.
Since $d\mid m$ and $d\mid\alpha$, we have $d\mid g$ and
$g\geq d>1$. By Lemma~\ref{lem:gcd-block},
\[
\{g-1,g,\ldots,2g-1\}\subseteq Z(m,\alpha)+\{0,1\}.
\]
Choose $s\in\{0,\ldots,d-1\}$ such that
$s\equiv r-g+1\pmod d$, and put $z=g-1+s$.
Then $g-1\leq z\leq2g-2$, and hence
\[
z\in Z(m,\alpha)+\{0,1\}
\subseteq\bigl(Z(m,\alpha)\cup Z(m,\alpha+1)\bigr)+\{0,1\}.
\]
Moreover, $z\equiv r\pmod d$, as required.

Now suppose that $\alpha=0$. We first show that
\[
\{0,2,\ldots,m-1\}\subseteq Z(m,1).
\]
For $i=0$, the permutation $\tau_0\lambda_1=\lambda_1$
is an $m$-cycle. Since
\[
x^{\tau_i\lambda_1}=
\begin{cases}
x+2\pmod m, & 0\leq x<i,\\
1, & x=i,\\
x+1\pmod m, & i<x<m,
\end{cases}
\]
for even $i$ with $0<i<m$, we have
\[
\tau_i\lambda_1
=(0,2,\ldots,i,1,3,\ldots,i-1,i+1,i+2,\ldots,m-1),
\]
where the sequence $i+1,i+2,\ldots,m-1$ is omitted when $i=m-1$. Thus $i\in Z(m,1)$.
If $r$ is even, then $r\in Z(m,1)$; if $r$ is odd,
then $r-1\in Z(m,1)$. In either case,
\[
r\in Z(m,1)+\{0,1\}
\subseteq\bigl(Z(m,0)\cup Z(m,1)\bigr)+\{0,1\}.
\]
Taking $z=r$ completes the proof.
\end{proof}

Throughout the rest of this section we use the notation of Lemma~\ref{lem:triple}. For $j\in\{1,2,3\}$, let
\begin{equation}\label{eq:window}
U_j=\bigcup_{\kappa=0}^{h_j-1}Z(m_j,\alpha_j+\kappa),
\qquad F_j=U_j\cup(U_j+1)\subseteq\Z_{m_j}.
\end{equation}
The next lemma gives the required residue-class coverage.

\begin{lemma}\label{thm:window}
For every ordered pair $i\ne j$, the set $F_j$ meets every residue class modulo $h_i$.
\end{lemma}

\begin{proof}
Fix distinct $i,j,k\in\{1,2,3\}$. Since $B_i+B_j=A$, the second isomorphism theorem gives
\[
B_j/(B_i\cap B_j)\cong A/B_i,
\]
and hence
\[
|B_j:B_i\cap B_j|=|A:B_i|=h_i.
\]
In particular, $h_i\mid m_j$, since $|B_j|=m_j$.

We next show that the assertion is unchanged when $s_i$ and $s_k$
are interchanged. Suppose that
\[
\delta_j=s_k-s_i,
\qquad
h_j(s_j-s_k)=\alpha_j\delta_j.
\]
After interchanging $s_i$ and $s_k$, we have
\[
\delta'_j=s_i-s_k=-\delta_j
\]
and
\[
h_j(s_j-s_i)
=h_j(s_j-s_k)+h_j(s_k-s_i)
=(\alpha_j+h_j)\delta_j
=(-\alpha_j-h_j)\delta'_j.
\]
Thus the new parameters are
\[
\delta'_j=-\delta_j,
\qquad
\alpha'_j=-\alpha_j-h_j\quad\text{in }\Z_{m_j}.
\] Let $U'_j$ and $F'_j$ be the corresponding sets in~\eqref{eq:window}.
By Lemma~\ref{lem:recurrence}, we have
\[
Z(m_j,-a-1)=-1-Z(m_j,a)
\qquad\text{for every }a\in\Z_{m_j}.
\]
Since $\alpha'_j=-\alpha_j-h_j$, setting
$\kappa'=h_j-1-\kappa$ gives
\[
\alpha'_j+\kappa'=-\alpha_j-\kappa-1.
\]
Consequently,
\[
\begin{aligned}
U'_j
&=\bigcup_{\kappa'=0}^{h_j-1}Z(m_j,\alpha'_j+\kappa')\\
&=\bigcup_{\kappa=0}^{h_j-1}Z(m_j,-\alpha_j-\kappa-1)\\
&=\bigcup_{\kappa=0}^{h_j-1}\bigl(-1-Z(m_j,\alpha_j+\kappa)\bigr)\\
&=-1-U_j.
\end{aligned}
\]
It follows that
\[
F'_j
=U'_j\cup(U'_j+1)
=(-1-U_j)\cup(-U_j)
=-\bigl(U_j\cup(U_j+1)\bigr)
=-F_j.
\]
Clearly, either both $F_j$ and $F'_j$ meet every residue class modulo $h_i$, or neither does.

We may therefore assume that $\delta_j=s_k-s_i$, in which case $h_j(s_j-s_k)=\alpha_j\delta_j.$
Since $B_i=\langle s_j-s_k\rangle$, we have $h_j(s_j-s_k)\in B_i.$
Moreover, $|A/B_j|=h_j$, so $
h_j\bigl((s_j-s_k)+B_j\bigr)=B_j,
$
and hence $h_j(s_j-s_k)\in B_j$. Therefore,
\[
\alpha_j\delta_j=h_j(s_j-s_k)\in B_i\cap B_j.
\]
Since $B_j=\langle\delta_j\rangle$ and $|B_j:B_i\cap B_j|=h_i$, it follows that $\alpha_j\equiv0\pmod{h_i}$. By Lemma~\ref{lem:triple}, $m_j$ is odd and $h_i,h_j\geq3$. Applying Lemma~\ref{lem:aligned-window} with $m=m_j$ and $d=h_i$, we conclude that
\[
\bigl(Z(m_j,\alpha_j)\cup Z(m_j,\alpha_j+1)\bigr)+\{0,1\}
\]
meets every residue class modulo $h_i$. This set is contained in $F_j$ by~\eqref{eq:window}, completing the proof.
\end{proof}

\begin{proof}[\rm\textbf{Proof of Proposition~\ref{prop:threefold}}]
Let $S=\{s_1,s_2,s_3\}$, and fix distinct $i,j\in\{1,2,3\}$. By Lemma~\ref{lem:triple}, there exists $Y_i\subseteq\mathrm{HP}(A,S)$ such that $2Y_i=2s_i+B_i$.
For each $\ell\in F_j$, the definition in~\eqref{eq:window}
gives $\ell=q+\varepsilon$ for some $0\leq\kappa<h_j$,
$q\in Z(m_j,\alpha_j+\kappa)$ and $\varepsilon\in\{0,1\}$.
By Lemma~\ref{lem:triple}, there is a Hamilton path from $0_+$
to $(s_j+(q+\varepsilon)\delta_j)_-$.
And so
\[
E_j=\{s_j+\ell\delta_j\mid\ell\in F_j\}
\subseteq\mathrm{HP}(A,S).
\]
Since $B_i+B_j=A$ and $B_j=\langle\delta_j\rangle$, the quotient $A/B_i$ is generated by $\delta_j+B_i$ and has order $h_i$. Lemma~\ref{thm:window} therefore implies that $E_j$ meets every coset of $B_i$. Consequently, $B_i+E_j=A$, and hence
\[
A=2s_i+B_i+E_j=2Y_i+E_j\subseteq3\mathrm{HP}(A,S)\subseteq A.
\]
This completes the proof.
\end{proof}

\section{Proof of main results}\label{sec:main-results}

Let $A$ be a finite abelian group and $T$ a non-empty subset of $A$. Recall that $D(T)$ denotes the subgroup $\langle s-t\mid s,t\in T\rangle$. We call the set $T$ a \emph{minimal $D$-set of $A$} if $D(T)=A$ but $D(T\setminus\{t\})\neq A$ for every $t\in T$. By Lemma~\ref{lem:haar}(d), the Haar graph $\Haar(A,T)$ is connected if and only if $D(T)=A$. If $D(T)\neq A$, then $\Haar(A,T)$ is disconnected and hence not Hamiltonian. If $D(T)=A$ but $T$ is not minimal, we may remove redundant elements from $T$ without affecting $D(T)$, and thus assume $T$ is minimal. Therefore, in studying the Hamiltonicity of Haar graphs over $A$, it suffices to consider minimal $D$-sets $T$ of $A$.

We first give some properties of minimal $D$-sets.

\begin{lemma}\label{lem:gdminimal}
Let $A$ be a finite abelian group, and let $T$ be a minimal $D$-set of $A$ with size at least $3$.  For each $t\in T$, let $L_t=D(T\setminus\{t\})$. Then the following hold.
\begin{enumerate}[label=\textup{(\alph*)},leftmargin=*]
\item\label{item:gdminimal-quotient} For any distinct $s,t\in T$, the quotient $A/L_t$ is cyclic of order greater than $1$ and generated by $t-s+L_t$.
\item\label{item:gdminimal-sum} $L_s+L_t=A$ for distinct $s,t\in T$.
\item\label{item:gdminimal-subsets} For each $I\subseteq T$ with $|I|\geq2$ and each $b\in I$, the set $I-b$ is a minimal $D$-set of $D(I)$.
\end{enumerate}
\end{lemma}

\begin{proof}

\ref{item:gdminimal-quotient}
Since $u-v=(u-s)-(v-s)$ for all $u,v\in T$, we have
\[
D(T)=\langle u-s\mid u\in T\rangle.
\]
Therefore, $$A=D(T)=\langle u-s\mid u\in T\setminus\{t\}\rangle +\langle t-s\rangle\leq L_t +\langle t-s\rangle.$$
Moreover, since $T$ is a minimal $D$-set of $A$, it follows that $A/L_t$ is cyclic
of order greater than $1$ and generated by $t-s+L_t$.

\ref{item:gdminimal-sum}
Since $|T|\geq 3$, we can choose $b\in T\setminus\{s,t\}$. Then, by \ref{item:gdminimal-quotient},
\[A\leq L_t+\langle t-b\rangle \leq L_t+L_s\leq A,\]
which proves the equality.


\ref{item:gdminimal-subsets}  Note that
\[
I-b\subseteq D(I),\qquad
D(I-b)=\langle (x-b)-(y-b)\mid x,y\in I\rangle=D(I),
\]
 Suppose that $D(I\setminus\{u\})=D(I)$ for some $u\in I$, and choose $v\in I\setminus\{u\}$. Then
\[
u-v\in D(I)=D(I\setminus\{u\})\subseteq D(T\setminus\{u\})=L_u.
\]
On the other hand, taking differences from $v$ gives
\[
A=D(T)=L_u+\langle u-v\rangle=L_u,
\]
contrary to the assumption. Thus, for every $u-b\in I-b$,  $D((I-b)\setminus\{u-b\})=D(I\setminus\{u\}) \neq D(I)$. This gives that $I-b$ is a minimal $D$-set of $D(I)$.
\end{proof}

\begin{theorem}\label{thm:endpoints}
Let $A$ be a finite abelian group of odd order, and let $T$ be a minimal $D$-set of $A$ with $|T|\geq4$.
Then
\[
2\mathrm{HP}(A,T)=A,
\]
and every edge of $\Haar(A,T)$ belongs to a Hamilton cycle.
\end{theorem}

\begin{proof}
We use induction on $r=|T|$. By Lemma~\ref{lem:haar}\ref{item:haar-components}, the assumptions give $D(T)=A$ and $D(T\setminus\{t\})\neq A$ for every $t\in T$.
Fix $t\in T$, choose $b\in T\setminus\{t\}$, and let
\[
L_t=D(T\setminus\{t\}),\qquad
U=(T\setminus\{t\})-b,\qquad s=t-b.
\]
Then $U\subseteq L_t$ and $D(U)=L_t$. By Lemma~\ref{lem:gdminimal}\ref{item:gdminimal-subsets},
$U$ is a minimal $D$-set of $L_t$.
Part~\ref{item:gdminimal-quotient} of the same lemma gives
\[
A=L_t+\langle s\rangle,\qquad m=|A:L_t|>1.
\]
Since $|A|$ is odd, $m$ is odd and hence $m\geq3$. Thus $A$, $L_t$, $s$ and $U$ satisfy Hypothesis~\ref{hyp:concat}.

If $r=4$, then $|U|=3$, and Proposition~\ref{prop:threefold} gives
\[
3\mathrm{HP}(L_t,U)=L_t.
\]
If $r\geq5$, then $|U|=r-1\geq4$, and the induction hypothesis gives
\[
2\mathrm{HP}(L_t,U)=L_t.
\]
In both cases, there is an integer $\ell\leq m$ such that $\ell\mathrm{HP}(L_t,U)=L_t$. Applying Lemma~\ref{lem:concat}\ref{item:concat-coverage}, we obtain
\[
s+L_t\subseteq\mathrm{HP}(A,U\cup\{s\})=\mathrm{HP}(A,T-b).
\]
By Lemma~\ref{lem:haar}\ref{item:haar-translation}, the isomorphism
\[
\Haar(A,T)\longrightarrow\Haar(A,T-b),\qquad
x_+\longmapsto x_+,\quad y_-\longmapsto(y-b)_-
\]
fixes $0_+$ and sends $a_-$ to $(a-b)_-$. Hence
\[
\mathrm{HP}(A,T-b)=\mathrm{HP}(A,T)-b.
\]
Adding $b$ to the inclusion above and using $s+b=t$, we conclude that
\begin{equation}\label{eq:deletion-endpoints}
t+L_t\subseteq\mathrm{HP}(A,T)\qquad(t\in T).
\end{equation}
For distinct $t,u\in T$, Lemma~\ref{lem:gdminimal}\ref{item:gdminimal-sum} now gives
\[
A=t+u+(L_t+L_u)
=(t+L_t)+(u+L_u)
\subseteq2\mathrm{HP}(A,T)\subseteq A.
\]
This proves $2\mathrm{HP}(A,T)=A$.

For each $t\in T$, we have $t\in t+L_t$, so~\eqref{eq:deletion-endpoints} gives a Hamilton path from $0_+$ to $t_-$. Since $|A|\geq|T|\geq4$, this path has more than two vertices and does not contain the edge $\{0_+,t_-\}$. Adding this edge gives a Hamilton cycle containing it. For each $a\in A$, translation by $a$ is an automorphism by Lemma~\ref{lem:haar}\ref{item:translation}, and sends this edge to
\[
\{a_+,(a+t)_-\}.
\]
Every edge of $\Haar(A,T)$ has this form for some $a\in A$ and $t\in T$. This proves the assertion about edges and completes the induction.
\end{proof}

\begin{proof}[\rm\textbf{Proof of Theorem~\ref{thm:generalized}}]
Let $\Gamma=\Cay(\Dih(A),S)$ be connected, where $|\Dih(A)|\geq4$ and $S\subseteq\Dih(A)\setminus\{1\}$. If $S\cap A\neq\varnothing$, then Pastine and Jaume~\cite[Theorem~15]{PJ2012} give a directed Hamilton cycle in $\Gamma$.

We may therefore assume that $S\subseteq\Dih(A)\setminus A$. Every element of $S$ is an involution, so $S=S^{-1}$ and $\Gamma$ is a Cayley graph. Since every edge can be traversed in either direction, it suffices to find a Hamilton cycle in this graph. If $|S|=2$, then $\Gamma$ is a cycle, so we may assume that $|S|\geq3$.

Suppose first that $|A|$ is even. Then $|\Dih(A)|=2|A|$ is divisible by $4$. By Alspach, Chen and Dean~\cite[Theorem~1.8]{ACD2010}, $\Gamma$ is Hamilton-connected or Hamilton-laceable. In either case, for any edge $\{u,v\}$ of $\Gamma$, there is
a Hamilton path $P$ from $u$ to $v$. Since $|V(\Gamma)|=|\Dih(A)|>2$, the edge $\{u,v\}$ does not belong
to $P$. Thus $P+\{u,v\}$ is a Hamilton cycle of $\Gamma$.

It remains to consider the case where $|A|$ is odd. Let
\[
S_0=\{a\in A\mid a\sigma\in S\}.
\]
By Lemma~\ref{lem:haar}\ref{item:haar-cayley}, $\Gamma\cong\Haar(A,S_0)$. Choose $T\subseteq S_0$ such that $\Haar(A,T)$ is connected
and $\Haar(A,T\setminus\{t\})$ is disconnected for every $t\in T$. Then $T$ is a minimal $D$-set of $A$.
Since $|A|>1$ and a Haar graph with at most one label is disconnected, $|T|\geq2$. If $|T|=2$, then $\Haar(A,T)$ is a connected graph of valency two and hence is a Hamilton cycle. If $|T|\geq4$, then Theorem~\ref{thm:endpoints} gives a Hamilton cycle in $\Haar(A,T)$.

Suppose that $|T|=3$. By Lemma~\ref{lem:haar}\ref{item:haar-translation}, we may assume that $T=\{0,t_1,t_2\}$. Let
\[
r=|t_1|,\qquad h=|A:\langle t_1\rangle|.
\]
We have $A=D(T)=\langle t_1,t_2\rangle$, so $t_2+\langle t_1\rangle$ generates $A/\langle t_1\rangle$ and has order $h$. In particular, $ht_2\in\langle t_1\rangle$. There is therefore a unique integer $\beta$ such that
\[
ht_2=\beta t_1,\qquad 0\leq\beta<r.
\]
By the minimal of $T$,  $r,h\geq3$ are odd.
Choose $0\leq\ell<2r$ such that
$\ell\equiv2\beta-h\pmod{2r}$. Then $h+\ell$ is even and
\[
\frac{h+\ell}{2}\equiv\beta\pmod r.
\]
Writing $A$ multiplicatively, let $x$ and $y$ correspond to
$-t_1$ and $-t_2$, respectively. Since $|A|=rh$, we have
\[
A\cong\langle x,y\mid x^r=1,\ xy=yx,\ y^h=x^\beta\rangle.
\]
Indeed, these relations reduce every element to the form
$x^iy^j$ with $0\leq i<r$ and $0\leq j<h$, and the natural
map onto $A$ is therefore an isomorphism. Consequently,
\[
\Dih(A)\cong
\left\langle x,y,\sigma\ \middle|\
\begin{gathered}
x^r=\sigma^2=1,\quad xy=yx,\quad y^h=x^{(h+\ell)/2},\\
\sigma x\sigma=x^{-1},\quad \sigma y\sigma=y^{-1}
\end{gathered}
\right\rangle.
\]
Under this identification, $T=\{0,t_1,t_2\}$ corresponds to
$\{1,x^{-1},y^{-1}\}$. Thus Lemma~\ref{lem:haar}\ref{item:haar-cayley}
and Alspach and Dean~\cite[Theorem~3.4]{AD2009} give
\[
\Haar(A,T)
\cong\Cay\bigl(\Dih(A),\{\sigma,x^{-1}\sigma,y^{-1}\sigma\}\bigr)
\cong\operatorname{HTG}(h,2r,\ell),
\]
where $\operatorname{HTG}(h,2r,\ell)$  is the honeycomb toroidal
graph with parameters $h$, $2r$ and $\ell$ defined in~\cite{AD2009}.
By~\cite[Theorem~4.1]{Alspach2021HTG}, $\operatorname{HTG}(h,2r,\ell)$ has a Hamilton cycle. It follows that $\Haar(A,T)$ has a Hamilton cycle.

In each case, $\Haar(A,T)$ has a Hamilton cycle. It is a spanning subgraph of $\Haar(A,S_0)$, so the same cycle is a Hamilton cycle of $\Haar(A,S_0)$. The isomorphism with $\Gamma$ gives a Hamilton cycle in the Cayley graph $\Gamma$, and hence a directed Hamilton cycle in the original digraph.
\end{proof}

\begin{remark}
Let $A$ be a finite abelian group of odd order, and let
$T\subseteq S\subseteq A$ with $|T|\geq4$. Suppose that
$T$ is a minimal $D$-set of $A$.
Since $\Haar(A,T)$ is a spanning subgraph of $\Haar(A,S)$,
we have
\[
\mathrm{HP}(A,T)\subseteq\mathrm{HP}(A,S).
\]
Theorem~\ref{thm:endpoints} therefore gives
\[
A=2\mathrm{HP}(A,T)\subseteq2\mathrm{HP}(A,S)\subseteq A,
\]
and hence $2\mathrm{HP}(A,S)=A$. Moreover, every edge of
$\Haar(A,T)$ belongs to a Hamilton cycle of $\Haar(A,S)$.
\end{remark}

\begin{proof}[\rm\textbf{Proof of Corollary~\ref{cor:dihedral}}]
Every dihedral group of order at least $4$ is a generalized dihedral group, so Theorem~\ref{thm:generalized} applies. The dihedral group of order $2$ gives a directed cycle of length $2$.
\end{proof}

\section*{AI Declaration}

ChatGPT was used to assist with mathematical exploration, proof development and revision, and the preparation of this manuscript. In earlier discussions with GPT-5.6, we investigated an approach to the Hamiltonicity problem for generalized dihedral Cayley graphs through Hamilton-laceability of honeycomb toroidal graphs. These explorations produced partial results but did not yield a complete solution. We later resumed the investigation with GPT-6, focusing first on the dihedral special case, namely the undirected case of Conjecture~\ref{conj:dihedral}. Building on the response from GPT-6, we obtained a proof that every connected Cayley graph on a dihedral group of order at least $4$ has a Hamilton cycle. This result prompted us to revisit the earlier approach and establish the corresponding result for generalized dihedral groups. Later, we provided GPT-6 with the idea of using the monodromy representation of graphs introduced by Yuan and Wang~\cite{YuanWang} to construct surface embeddings with a face bounded by a Hamilton cycle, which led to an essential simplification of the original proof. Finally, combining this undirected result with the theorem of Pastine and Jaume~\cite[Theorem~15]{PJ2012} gave the directed version in Theorem~\ref{thm:generalized}.

The authors directed the research and are responsible for the mathematical content and presentation.

\end{document}